\documentclass[twoside,notitlepage,12pt]{article}

\usepackage{amssymb}
\usepackage[leqno]{amsmath}
\usepackage{amsfonts}
\usepackage{amsopn}
\usepackage{amstext}
\usepackage{amsthm}
\usepackage{hyperref}
\usepackage{makeidx}

\newcommand{\Ef}{{\cal{F}}}

\newcommand{\Tau}{{\cal{T}}}

\newcommand{\Qyu}{{\Bbb{Q}}}
\newcommand{\Err}{{\Bbb{R}}}

\newcommand{\lam}{{\lambda}}
\newcommand{\al}{\alpha}
\newcommand{\Gam}{\Gamma}

\newcommand{\sig}{\sigma}
\newcommand{\eps}{\varepsilon}
\renewcommand{\phi}{\varphi}
\renewcommand{\rho}{\varrho}

\newcommand{\rest}{\restriction}

\newcommand{\ntr}{n\in\omega}

\newcommand{\loe}{\leqslant}
\newcommand{\goe}{\geqslant}

\newcommand{\subs}{\subseteq}
\newcommand{\sups}{\supseteq}
\newcommand{\nnempty}{\ne\emptyset}
\newcommand{\argum}{\:\cdot\:}
\newcommand{\ovr}{\overline}
\renewcommand{\iff}{\Longleftrightarrow}

\newcommand{\cl}{\operatorname{cl}}

\newcommand{\id}[1]{\operatorname{id}_{#1}}
\newcommand{\cf}{\operatorname{cf}}

\newcommand{\pr}{\operatorname{pr}}

\newcommand{\by}{/}

\newtheorem{tw}{Theorem}[section]
\newtheorem{wn}[tw]{Corollary}
\newtheorem{lm}[tw]{Lemma}
\newtheorem{prop}[tw]{Proposition}

\theoremstyle{definition}

\newtheorem{ex}[tw]{Example}
\theoremstyle{remark}

\newcommand{\setof}[2]{\{#1\colon #2\}}

\newcommand{\sett}[2]{\{#1\}_{#2}}
\newcommand{\sn}[1]{\{#1\}} % singleton
\newcommand{\pair}[2]{\langle #1, #2 \rangle} % pair
\newcommand{\triple}[3]{\langle #1, #2, #3 \rangle} % triple
\newcommand{\map}[3]{#1\colon #2 \to #3} % A function
\newcommand{\img}[2]{#1[#2]} % image of a set
\newcommand{\inv}[2]{{#1}^{-1}[#2]} % preimage of a set

\newcommand{\dpower}[2]{[#1]^{#2}}

\providecommand{\cal}{\mathcal}
\renewcommand{\Bbb}{\mathbb}

\newcommand{\cmp}{\circ} % composition!!!
\newcommand{\s}{\mathfrak s}
\newcommand{\ft}{\mathfrak t}

\newenvironment{pf}{\begin{proof}}{\end{proof}}

\newcommand{\rloe}{\preceq}

\newcommand{\suppt}{\operatorname{suppt}}

\newcommand{\Sig}{\Sigma}

\providecommand{\nat}{\omega}

\newcommand{\im}{\operatorname{im}}
\newcommand{\bal}{\operatorname{B}}
\newcommand{\clbal}{\overline{\bal}}
\newcommand{\anorm}{\|\cdot\|}
\newcommand{\norm}[1]{\|#1\|}

\newcommand{\aabs}{|\cdot|}
\newcommand{\abs}[1]{|#1|}
\newcommand{\uball}[1]{\clbal_{#1}}
\newcommand{\ubal}{\uball}
\newcommand{\usphere}[1]{\operatorname{S}_{#1}}
\newcommand{\dens}{\operatorname{dens}}

\newcommand{\weakstar}{\ensuremath{{weak}^*}}

\newcommand{\R}{\ensuremath{\mathcal R}}
\newcommand{\Rzero}{\ensuremath{\mathcal R_0}}

\newcommand{\invsys}[5]{\langle {#1}_{#4},{#2}_{#4}^{#5},#3 \rangle}

\newcommand{\leftorth}[1]{{}^\perp{\left(#1\right)}}
\newcommand{\rightorth}[1]{{\left(#1\right)}^\perp}
\newcommand{\cee}[1]{{\mathcal C}\!\left({#1}\right)}
\newcommand{\ceep}[1]{{\mathcal C}_p\!\left({#1}\right)}

\newcommand{\clstar}{\cl_*}

\newcommand{\tc}[1]{\operatorname{tc}(#1)} % transitive closure
\newcommand{\tcn}[2]{\operatorname{tc}_{#2}(#1)} % transitive n-closure

\newtheorem{question}{Question}
\newtheorem{problem}{Problem}

\title{Banach spaces with projectional skeletons}
\author{
{\sc Wies{\l}aw Kubi\'s}
\footnote{Research supported by MNiSW grant No. N201 024 32/0904. Author's e-mail address: \texttt{wkubis@pu.kielce.pl}}
\\
{\small Institute of Mathematics, Jan Kochanowski University in Kielce, Poland}\\
{\it and}\\
{\small Mathematical Institute of the Czech Academy of Sciences, Prague, Czech Republic}
}

\makeindex

\begin{document}
\maketitle

\begin{abstract}
A projectional skeleton in a Banach space is a $\sig$-directed family of projections onto separable subspaces, covering the entire space. The class of Banach spaces with projectional skeletons is strictly larger than the class of Plichko spaces (i.e. Banach spaces with a countably norming Markushevich basis). We show that every space with a projectional skeleton has a projectional resolution of the identity and has a norming space with similar properties to $\Sig$-spaces. We characterize the existence of a projectional skeleton in terms of elementary substructures, providing simple proofs of known results concerning weakly compactly generated spaces and Plichko spaces.

We prove a preservation result for Plichko Banach spaces, involving transfinite sequences of projections. As a corollary, we show that a Banach space is Plichko if and only if it has a commutative projectional skeleton.


\end{abstract}

\newpage
\tableofcontents

\section{Introduction}

It is well known that a Banach space in which every closed subspace is complemented is necessarily isomorphic to a Hilbert space (Lindenstrauss \& Tzafriri \cite{LinTza71}). On the other hand, there exist Banach spaces in which only finite-dimensional and co-finite-dimensional subspaces are complemented (Gowers \& Maurey \cite{GowMau93}). There even exist (necessarily, non-separable) $\cee K$ spaces with no nontrivial bounded projections (Koszmider \cite{Kos04} and Plebanek \cite{Ple04}). In the positive direction, one has to mention the work of Heinrich \& Mankiewicz \cite{HeiMan} where, using substructures of ultrapowers of Banach spaces, the authors show in particular the existence of non-trivial bounded projections in every dual Banach space of density greater than the continuum (see \cite{SimYos} for an elementary proof).

We are interested in Banach spaces which have ``many" projections onto separable subspaces. Perhaps one of the most general classes of this sort would be the class of Banach spaces with the {\em separable complementation property} (SCP): a space $X$ has the SCP if, by definition, for every countable set $S\subs X$ there is a bounded projection $\map PXX$ such that $\im P$ is separable and $S\subs \im P$. It turns out that this property is general enough in order to include somewhat pathological spaces.
For a survey on the complementation property and its variations we refer to \cite{PlichkoYost00}.

Another possible class of Banach spaces with many projections are spaces with the {\em projectional resolution of the identity} (PRI), notion introduced by Lindenstrauss \cite{Lin65, Lin66}, defined to be a well ordered continuous chain of projections onto smaller subspaces. This property together with transfinite induction allows proving various properties of a non-separable Banach space, e.g. a locally uniformly rotund/Kadec renorming \cite{Ziz84, BurKubTod} and the existence of a linear injection into $c_0(\Gam)$ \cite{AmiLin68, Vas81}. See \cite{Fab97} or \cite{DevGodZiz} for more information concerning the PRI method.
Unfortunately, the existence of a PRI in a Banach space is good enough only for density $\aleph_1$, otherwise it does not even imply the SCP.

We propose a natural class of Banach spaces which have a family of projections onto separable subspaces, indexed by a $\sig$-directed partially ordered set and satisfying some natural conditions, similar to a PRI. We call it a ``projectional skeleton". It turns out that a Banach space with a projectional skeleton looks ``almost" like a Plichko space, i.e. a Banach space with a countably norming Markushevich basis (see \cite{Pli82, Pli83, Pli86} and \cite{Val91}). In fact, we essentially know only one basic example distinguishing those two classes: the space $\cee K$, where $K$ is the ordinal $\omega_2+1$, endowed with the interval topology. This was shown by Kalenda in \cite{Kal02}. Banach spaces with a projectional skeleton can be characterized by a property involving norming sets and countable elementary submodels. We explain
%in Section \ref{wejfopaqwi}
how to use elementary submodels of (some initial part of) set theory for constructing bounded projections and we show that the existence of a projectional skeleton is equivalent to some natural model-theoretic property of a suitable norming space, similar to the existence of the so-called {\it projectional generator}. As an application, we get short proofs of some well known results, like the existence of a PRI in every weakly compactly generated space. Our characterization allows us to show that the class of Banach spaces with a projectional skeleton of norm one is stable under arbitrary $c_0$- and $\ell_p$-sums ($1\loe p<\infty$). We apply elementary submodels for constructing a PRI from a projectional skeleton.

The main part (Section \ref{popop}) is devoted to Plichko spaces. We prove a preservation result for inductive limits of certain projective sequences of Plichko spaces, similar in spirit to Gul$'$ko's results on subspaces of the $\Sig$-products \cite{Gul77, Gul79, Gul98}. As an application, we show that a Banach space is Plichko if and only if it has a commutative projectional skeleton.

Finally, we discuss retractional skeletons in compact spaces -- a notion dual to projectional skeleton. We characterize this class of compacta by means of elementary submodels and we state a preservation property for Valdivia compacta, dual to the corresponding result for Banach spaces. Retractional skeletons were introduced already in \cite{KubMich}, where it is proved that Valdivia compacta are precisely those compact spaces which have a commutative retractional skeleton. For more information and recent results concerning Valdivia compacta and their spaces of continuous functions we refer to \cite{Kal00a, Kal-natural, KubMich, KubUsp, BanKub, Kub06}.

\section{Preliminaries}

We shall consider Banach spaces over the field of real numbers, although the results are true also for the complex case. Below we recall most relevant notions, definitions and notation.

\index{projection}
By a {\em projection} in a Banach space $X$ we mean a bounded linear operator $\map PXX$ such that $P\cmp P=P$. In this case $\im P = \setof{x\in X}{x = Px}$ and $\ker P = \setof{x-Px}{x\in X} = \im(\id X - P)$, where $\id X$ is the identity map.
\index{separable complementation property}\index{SCP}
Recall that a space $X$ has the {\em separable complementation property} ({\em SCP} for short) if for every countable set $S\subs X$ there is a projection $\map PXX$ such that $PX$ is a separable space containing $S$.
Given $B\subs X^*$, we write $\leftorth B = \setof{x\in X}{(\forall\;b\in B)\;b(x)=0}$. The right annihilator $\rightorth A$ is defined similarly.

\index{projectional sequence}
Let $\lam$ be a limit ordinal. A {\em projectional sequence} of length $\lam$ in a Banach space $X$ is a sequence of projections $\sett{P_\xi}{\xi<\lam}$ satisfying the following conditions:
\begin{enumerate}
	\item $\xi<\eta\implies P_\xi=P_\eta\cmp P_\xi=P_\xi\cmp P_\eta$,
	\item $P_\delta X=\cl(\bigcup_{\xi<\delta}P_\xi X)$ for every limit ordinal $\delta<\lam$,
	\item $X=\cl(\bigcup_{\xi<\lam}P_\xi X)$.
\end{enumerate}
\index{projectional resolution of the identity}
A special case is a {\em projectional resolution of the identity} ({\em PRI}): this is a projectional sequence $\sett{P_\xi}{\xi<\lam}$ such that $\norm{P_\xi}=1$ and $\dens(P_\xi X)\loe|\xi|+\aleph_0$ for each $\xi<\lam$, where $|\xi|$ denotes the cardinality of $\xi$.

%\paragraph{Topology.}
All topological spaces are assumed to be completely regular.
\index{space}
\index{countably tight space}\index{space!-- countably tight}
The closure of a set $A$ in a space $X$ will be denoted by $\cl(A)$ or, more precisely, by $\cl_X(A)$. If $X$ is a dual to a Banach space then $\clstar$ will denote the \weakstar\ closure, i.e. the closure with respect to the \weakstar\ topology on $X$.

A space $X$ is {\em countably tight} if for every $A\subs X$ and for every $p\in \cl A$ there exists $A_0\in\dpower A{\aleph_0}$ with $p\in\cl A_0$.
\index{$\Sig$-product}
Let $\Gam$ be a set. Given $x\in \Err^\Gam$, we denote by $\suppt(x)$ the {\em support} of $x$, i.e. $\suppt(x)=\setof{\al\in\Gam}{x(\al)\ne0}$. The set $\Sig(\Gam) = \setof{x\in\Err^\Gam}{|\suppt(x)|\loe\aleph_0}$ is called a {\em $\Sig$-product}.

\index{Valdivia compact}\index{Corson compact}\index{$\Sig$-subset}
A {\em Valdivia compact} is a compact space homeomorphic to $K\subs[0,1]^\kappa$ satisfying $K=\cl(K\cap \Sig(\kappa))$. A {\em Corson compact} is, by definition, a compact subset of $\Sig(\kappa)$ for some $\kappa$. Given a compact $K$, $D\subs K$ is called a {\em $\Sig$-subset} of $K$ if there is a homeomorphic embedding $\map hK{[0,1]^\kappa}$ such that $D = \inv h{\Sig(\kappa)}$.

%\paragraph{Norming sets.}
\index{norming space/set}\index{space!-- norming}
Let $\pair X\anorm$ be a Banach space. A set $D\subs X^*$ is {\em norming} if the formula
\begin{equation}
\abs x = \sup\setof{|\phi(x)|/\norm\phi}{\phi\in D\setminus\sn0}
\tag{*}\label{gwiazda1}\end{equation}
defines an equivalent norm on $X$. More precisely, we say that $D$ is {\em $r$-norming} if $\norm x\loe r\abs x$ for every $x\in X$.
$D$ is {\em $1$-norming} if $\aabs = \anorm$.
The following fact is well known.

\begin{prop}\label{nofjsopfj}
Let $D$ be a norming subset of $X^*$. Then $D$ is $1$-norming with respect to the norm $\aabs$ defined by equation (\ref{gwiazda1}).
\end{prop}

\begin{pf}
Let $D'$ be the linear span of $D$ and let $D_1=\setof{\phi\in D'}{\norm\phi\loe1}$.
Then $\abs x = \sup_{\phi\in D_1}|\phi(x)|$. Thus, it remains to notice that $D_1=\setof{\phi\in D'}{\abs \phi\loe1}$. Indeed, if $\phi \in D_1$ then $|\phi(x)|\loe1$ whenever $\abs x\loe1$, so $\abs\phi\loe1$. Conversely, if $\phi\in D\setminus D_1$ then there is $x_0\in X$ such that $\abs {x_0}=1$ and $|\phi(x_0)|>1$, so $\abs\phi>1$.
\end{pf}

%\paragraph{Plichko spaces.}
\index{Plichko space}\index{Banach space!-- Plichko}
Recall that a Banach space $X$ is called {\em Plichko} if there exists a one-to-one \weakstar\ continuous linear operator $\map T{X^*}{\Err^\kappa}$ such that $\inv T{\Sig(\kappa)}$ is norming for $X$.
Equivalently: there are a linearly dense set $A\subs X$ and a norming set $D\subs X^*$ such that for every $y\in D$ the set $\setof{a\in A}{y(a)\ne0}$ is countable. If additionally $D$ is linear, we shall say that $\pair XD$ is a {\em Plichko pair}.
More generally, we say that $\pair YD$ is a Plichko pair in a space $X$ if $Y$ is a closed linear subspace of $X$ and $\pair Y{D'}$ is a Plichko pair, where $D'=\setof{y\rest Y}{y\in D}$.
Given $A\subs X$, the set $\suppt(y,A) = \setof{a\in A}{y(a)\ne0}$ will be called the {\em $A$-support} of $y\in X^*$. The space $\ovr D=\setof{y\in X^*}{|\suppt(y,A)|\loe\aleph_0}$ is called a {\em $\Sig$-space}.
\index{Plichko pair}\index{support}

\index{Banach space!-- weakly Lindel\"of determined}
\index{Banach space!-- WLD}
\index{WLD space}
A particularly interesting subclass of Plichko spaces is the class of weakly Lindel\"of determined spaces, introduced by Valdivia \cite{Val88}. A Banach space $X$ is {\em weakly Lindel\"of determined} ({\em WLD} for short) if $\pair X{X^*}$ is a Plichko pair, i.e. $X^*$ is a $\Sig$-space. This is equivalent to saying that the dual unit ball with the \weakstar\ topology is Corson compact (see \cite[Prop. 4.1]{OriSchVal91}).

\section{Elementary submodels and projections}\label{wejfopaqwi}

In this section we introduce the method of elementary submodels, which will be used extensively throughout the paper. In the context of retractions -- topological counterparts of linear projections -- elementary submodels turned out to be an important tool in \cite{Kub06, KubMich}. We refer to the survey article of Dow \cite{Dow88}, where several applications of elementary submodels in set-theoretic topology are explained.

Let $N$ be a fixed set. The pair $\pair N\in$, where $\in$ is restricted to $N\times N$, is a structure of the language of set theory. Given a formula $\phi(x_1,\dots,x_n)$ with all free variables shown and given $a_1,\dots,a_n\in N$ one defines the relation ``$\pair N\in$ satisfies $\phi(a_1,\dots,a_n)$" (briefly ``$\pair N\in\models\phi(a_1,\dots,a_n)$" or just ``$N\models \phi(a_1,\dots,a_n)$") in the usual way, by induction on the length of the formula. Namely, $N\models a_1\in a_2$ iff $a_1\in a_2$ and $N\models a_1=a_2$ iff $a_1=a_2$. It is clear how ``satisfaction" is defined for conjunction, disjunction and negation. Finally, if $\phi$ is of the form $(\exists\; y)\psi(x_1,\dots,x_n,y)$ then $N\models \phi(a_1,\dots,a_n)$ iff there exists $b\in N$ such that $N\models \psi(a_1,\dots,a_n,b)$. 

As an example, if $s=\{a,b,c\}$ and $s,a,b\in N$ while $c\notin N$, then $N$ satisfies ``$s$ has at most two elements", because for every $x\in N$ if $x\in s$ then either $x=a$ or $x=b$.

\index{relativization}
Instead of the above definition, some authors use relativization, see e.g. Kunen's book \cite{Kunen-book}. Given a formula $\phi$, the {\em relativization} of $\phi$ to $N$ is a formula $\phi^N$ which is built from $\phi$ by replacing each quantifier of the form ``$\forall\;x$" by ``$\forall\; x\in N$" and each quantifier of the form ``$\exists\; x$" by ``$\exists\;x\in N$". By this way, $N\models\phi(a_1,\dots,a_n)$ iff $\phi^N(a_1,\dots,a_n)$ holds (of course, $a_1,\dots,a_n$ must be elements of $N$).

\index{transitive closure}
\index{transitivity}
Given a set $x$, we define the {\em transitive closure} of $x$ to be $\tc x = \bigcup_{\ntr}\tcn xn$, where $\tcn x1= x\cup\bigcup x$ and $\tcn x{n+1} = \tcn{\tcn xn}1$. In other words: $y\in\tc x$ iff there are $x_0\in x_1\in \dots\in x_k$ such that $y\in x_0$ and $x_k\in x$. Thanks to the Axiom of Regularity, these two definitions of transitive closure are equivalent and every set of the form $\tc x$ is {\em transitive}, i.e. $y\in \tc x$ implies $y\subs \tc x$.

Given a cardinal $\theta$, we denote by $H(\theta)$ the class of all sets whose transitive closure has cardinality $<\theta$. It is well known that $H(\theta)$ is a set, not a proper class. It is clear that $H(\theta)$ is transitive. We shall consider elementary substructures of $\pair{H(\theta)}\in$. It is well known that for a regular uncountable cardinal $\theta$, the structure $\pair{H(\theta)}\in$ satisfies all the axioms of set theory, except possibly the Power Set Axiom, see \cite[IV.3]{Kunen-book}.

\index{elementary submodel/substructure}
Recall that a substructure $M$ of $\pair{H(\theta)}\in$ is called {\em elementary} if for every formula $\phi(x_1,\dots,x_n)$ with all free variables shown, for every $a_1,\dots, a_n\in M$, we have that
$$M\models \phi(a_1,\dots,a_n) \iff H(\theta)\models \phi(a_1,\dots,a_n).$$
The fact that $M$ is an elementary submodel of $\pair{H(\theta)}\in$ is denoted by $M\rloe \pair{H(\theta)}\in$.

%The sets $H(\theta)$ provide a hierarchy of the universe of set theory. 

The reason for using elementary submodels of $H(\theta)$ is that these structures 
satisfy most of the axioms of set theory: if $\theta>\aleph_0$ is regular then $H(\theta)$ satisfies all the axioms except possibly the power-set, because it may happen that $2^\lam>\theta$ for some $\lam<\theta$. Moreover, in practice it is usually easy to point out a cardinal $\theta$ such that $H(\theta)$ satisfies given finitely many formulas with parameters, needed for applications. Another useful feature of $H(\theta)$ with $\theta$ regular is the fact for every formula $\phi(x_1,\dots,x_n)$ in which all quantifiers are bounded (i.e. of the form ``$\forall\;x\in y$" or ``$\exists\;x\in y$") and for every $a_1,\dots,a_n\in H(\theta)$, $\phi(a_1,\dots,a_n)$ holds if and only if $H(\chi)\models \phi(a_1,\dots,a_n)$. For more information, see  \cite[IV.3]{Kunen-book}. Since in most cases we indeed use formulas with bounded quantifiers, one can simply ``check" their validity by looking at a sufficiently large $H(\theta)$.

One can also use
{\em Reflection Principle}, which says that given a formula of set theory $\phi(x_1,\dots,x_n)$ and given sets $a_1,\dots,a_n$ such that $\phi(a_1,\dots,a_n)$ holds, there exists $\theta$ such that the structure $\pair{H(\theta)}{\in}$ satisfies $\phi(a_1,\dots,a_n)$.
In some cases $\theta$ may not be regular, although it may be arbitrarily big and it may have arbitrarily big cofinality. More precisely: the class of cardinals $\theta$ with the above property is closed and unbounded.
%In fact, Reflection Principle says about a closed and unbounded class of cardinals, not just about one cardinal.
Thus, when considering finitely many formulas and parameters, we can ``check" their validity by restricting attention to $H(\theta)$, where $\theta$ is a ``big enough" cardinal, meaning that the cofinality of $\theta$ is greater than a prescribed cardinal and all relevant formulas are satisfied in $\pair{H(\theta)}\in$.
\index{$H(\theta)$}
\index{Reflection Principle}

Summarizing: assume we would like to use in our arguments formulas $\phi_1,\dots,\phi_n$ and parameters from a finite set $S$. We then find a cardinal $\theta$ so that $S\subs H(\theta)$ and, by Reflection Principle, all valid formulas $\phi_1,\dots,\phi_n$ with suitable parameters are satisfied in $H(\theta)$. Finally, we shall use elementary substructures of $H(\theta)$ which contain $S$. If the formulas $\phi_1,\dots,\phi_n$ have only bounded quantifiers, then we do not need to use Reflection Principle, since the formulas will be satisfied in every $H(\theta)$ with $\theta$ big enough, i.e. every regular $\theta$ greater than some fixed cardinal $\theta_0$.

In order to illustrate elementarity, let us come back to the simple example described above: let $s=\{a,b,c\} \in N$, $a,b\in N$ and now assume that $N\rloe H(\theta)$ and that $a,b,c$ are pairwise distinct. 
Since $N\subs H(\theta)$, we see that $s\in H(\theta)$ and consequently also $c\in H(\theta)$.
Clearly, $H(\theta)\models c\in s$, therefore $H(\theta)$ satisfies ``there is $x\in s$ such that neither $x=a$ nor $x=b$ ". By elementarity, $N$ satisfies the same statement, which means that there exists $d\in N$ such that $N\models (d\in s\land d\ne a\land d\ne b)$. This is a conjunction of atomic formulas and their negations, so indeed $d\in s$ and $d\notin\{a,b\}$. But we have assumed that $s\cap N=\{a,b\}$, which is a contradiction. This example shows that elementary substructures of $H(\theta)$ ``keep" elements of a finite set. In general, if $N\rloe H(\theta)$ then $s\in N$ does not necessarily imply that $s\subs N$, unless $s$ is countable or $N$ contains a sufficiently big initial interval of ordinals, see Proposition \ref{enbgojhe}(c) below or \cite[Thm. 1.6]{Dow88}.

A particular case of the L\"owenheim-Skolem Theorem (for the language of set theory) says that for every infinite set $S\subs H(\theta)$ there exists $M\rloe \pair{H(\theta)}\in$ such that $|M|=|S|$. This theorem can be viewed as the ``ultimate" closing-off argument and its typical proof indeed proceeds by ``closing-off" the given set $S$, by adding elements which witness ``satisfaction" of all suitable formulas of the form $(\exists\; x)\;\psi$.

Important for applications is the fact that, thanks to the L\"owenheim-Skolem theorem, we may consider \emph{countable} elementary substructures of an arbitrarily large $H(\theta)$.

\begin{prop}\label{enbgojhe}
Let $\theta$ be an uncountable regular cardinal and let $M\rloe \pair{H(\theta)}\in$.
\begin{enumerate}
	\item[(a)] Assume $u\in H(\theta)$, $a_1,\dots,a_n\in M$ and $\phi(y,x_1,\dots,x_n)$ is a formula such that $u$ is the unique element of $H(\theta)$ for which $H(\theta)\models \phi(u,a_1,\dots,a_n)$. Then $u\in M$.
	\item[(b)] Let $s\subs M$ be a finite set. Then $s\in M$.
	\item[(c)] If $S\in M$ is a countable set then $S\subs M$.
\end{enumerate}
\end{prop}

\begin{pf}
(a) By elementarity there exists $v\in M$ such that $M\models \phi(v,a_1,\dots,a_n)$. Using elementarity again, we see that $H(\theta)\models \phi(v,a_1,\dots,a_n)$. Thus $u=v$.

(b) Let $s=\{a_1,\dots,a_n\}$. Then $s$ is the unique set satisfying the formula $\phi(s,a_1,\dots,a_n)$, where $\phi(x,y_1,\dots,y_n)$ is
$$(\forall\;t)\; t\in x \iff t=y_1\lor t=y_2\lor \dots \lor t=y_n.$$
Applying (a), we see that $s\in M$.

(c) By induction and by (a), we see that all natural numbers are in $M$.
Also by (a), the set of natural numbers $\omega$ is an element of $M$, being uniquely defined as the minimal infinite ordinal. Notice that $H(\theta)$ satisfies ``there exists a surjection from $\omega$ onto $S$\/". By elementarity, there exists $f\in M$ such that $M$ satisfies ``$f$ is a surjection from $\omega$ onto $S$\/". Again using (a), we see that $f(n)\in M$ for each $\ntr$. Finally, it suffices to observe that $f$ is indeed a surjection, i.e. for every $x\in S$ there is $n$ such that $x=f(n)$. This follows from elementarity, because assuming $\img f\omega\ne S$, the formula ``$(\exists\;x\in S) (\forall\; \ntr)\;\;x\ne f(n)$" would be satisfied in $M$, contradicting that $f$ is a surjection.
\end{pf}

Fix a Banach space $X$ and choose a ``big enough" regular cardinal $\theta$, so that $X\in H(\theta)$. Take an elementary substructure $M$ of $\pair{H(\theta)}\in$ such that $X\in M$. Note that $M$ may be countable, by the L\"owenheim-Skolem Theorem. What can we say about the set $X\cap M$? By elementarity, it is closed under addition.
By Proposition \ref{enbgojhe}(a), the field of rationals is contained in $M$, therefore $X\cap M$ is a $\Qyu$-linear subspace of $X$. Consequently, the norm closure of $X\cap M$ is a Banach subspace of $X$. In particular, the weak closure of $X\cap M$ equals the norm closure of $X\cap M$.
We shall write $X_M$ instead of $\cl(X\cap M)$ and we shall call $X_M$ the {\em subspace induced by} $M$.

In case of some typical Banach spaces, we can describe the subspace $X_M$. For instance, let $X=\ell_p(\Gamma)$, where $1\loe p<\infty$ and $\Gamma$ is an uncountable set. Then $X_M$ can be identified with $\ell_p(\Gamma\cap M)$. Indeed, identify $x\in \ell_p(\Gamma\cap M)$ with its extension $x'\in\ell_p(\Gamma)$ defined by $x'(\al)=0$ for $\al\in \Gamma\setminus M$. Let $x\in X\cap M$. Then $\suppt(x)=\setof{\al\in\Gamma}{x(\al)\ne0}$ is a countable set and hence, by elementarity, it belongs to $M$. By Proposition \ref{enbgojhe}(c), $\suppt(x)\subs M$. Thus $x\in\ell_p(\Gamma\cap M)$. On the other hand, if $x\in \ell_p(\Gamma\cap M)$ then arbitrarily close to $x$ we can find $y\in \ell_p(\Gamma\cap M)$ such that $s=\suppt(y)$ is finite. Moreover, we may assume that $y(\al)\in\Qyu$ for $\al\in s$. By Proposition \ref{enbgojhe}(b), $y\rest s\in M$ and consequently also $y\in M$. Hence $x\in \cl (X\cap M)=X_M$.

Given a compact space $K\in H(\theta)$ and $M\rloe \pair{H(\theta)}\in$, define the following equivalence relation $\sim_M$ on $K$:
$$x\sim_M y\iff (\forall\;f\in\cee K \cap M)\;f(x)=f(y).$$
We shall write $K\by M$ instead of $K\by {\sim_M}$ and we shall denote by $q^M$ (or, more precisely, $q^M_K$) the canonical quotient map.
It is not hard to check that $K\by M$ is a compact Hausdorff space of weight not exceeding the cardinality of $M$. This construction has been used by Bandlow \cite{Ban91, Ban94} for characterizing Corson compacta in terms of elementary substructures.

\begin{lm}\label{ergnosfoipq3r}
Let $K$ be a compact space, let $\theta$ be a big enough regular cardinal and let $M\rloe \pair{H(\theta)}\in$ be such that $K\in M$. Then
$$\cl(\cee K\cap M)=\setof{\phi\cmp q^M}{\phi\in \cee {K\by M}},$$
where $\cl$ denotes the norm closure in the above formula.
\end{lm}

\begin{pf}
Let $Y$ denote the set on the right-hand side. Then $Y$ is a closed linear subspace of $\cee K$. Given $\psi\in\cee K\cap M$, by the definition of $\sim_M$, there exists a (necessarily continuous) function $\psi'$ such that $\psi=\psi'\cmp q^M$. Thus $\cee K\cap M\subs Y$. Let $R=\setof{\phi\in\cee{K\by M}}{\phi\cmp q^M\in M}$. Then $R$ is a subring of $\cee{K\by M}$ which separates points and which contains all rational constants. By the Stone-Weierstrass Theorem, $R$ is dense in $\cee{K\by M}$, which implies that $\cee K\cap M$ is dense in $Y$.
\end{pf}

Observe that, under the assumptions of the above Lemma, the norm closure of $\cee K\cap M$ is pointwise closed. Indeed, if $f\in \cee K\setminus \cl(K\cap M)$ then there are $x,y \in K$ such that $x\sim_M y$ while $f(x)\ne f(y)$. Consequently, $V=\setof{g}{g(x)\ne g(y)}$ is a neighborhood of $f$ in the pointwise convergence topology, disjoint from $\cl(K\cap M)$.

\subsection{Projections induced by elementary substructures}

Next we show how to use elementary submodels for constructing bounded projections.
This idea has already been applied, in case of WCG spaces, by Koszmider \cite{Kosz-wcg}.

\begin{lm}\label{sgagaqqwfwr}
Assume $X$ is a Banach space, $D\subs X^*$ is $r$-norming and $M$ is an elementary substructure of a big enough $\pair{H(\theta)}{\in}$ such that $X,D\in M$. Then 
\begin{enumerate}
	\item[(a)] $X_M\cap \leftorth{D\cap M}=\sn0$;
	\item[(b)] the canonical projection $\map P{X_M\oplus \leftorth{D\cap M}}{X_M}$ has norm $\loe r$.
\end{enumerate}
\end{lm}

\begin{pf}
Fix $x\in X\cap M$, $y\in \leftorth{D\cap M}$ and fix $\eps>0$.
Since $D$ is $r$-norming, there exists $d\in D$ such that $r|d(x)|/\norm d\goe \norm x-\eps$.
Since $x\in M$, by elementarity we may assume that $d\in M$. Thus $d\in D\cap M$ and $d(y)=0$. It follows that
$$\norm x\loe r|d(x)|/\norm d+\eps= r|d(x+y)|/\norm d+\eps\loe r\norm{x+y}+\eps.$$
By continuity, we see that $\norm x\loe r\norm{x+y}$ whenever $x\in X_M$ and $y\in\leftorth{D\cap M}$. In particular, $X_M\cap \leftorth{D\cap M}=\sn0$, because if $x\in X_M\cap \leftorth{D\cap M}$ then $-x\in \leftorth{D\cap M}$ and $\norm x\loe r\norm{x-x}=0$.
\end{pf}

Note that, in the above lemma, the subspace $X_M\oplus\leftorth{D\cap M}$ is closed in $X$.

It may happen that $\leftorth{D\cap M}=0$ (consider $X=\ell_\infty$) and in that case the above lemma is meaningless.
We are going to discuss Banach spaces for which Lemma \ref{sgagaqqwfwr} provides a way to construct full projections.

\subsection{WCG spaces and Plichko pairs}

We demonstrate the use of elementary submodels for finding projections in weakly compactly generated spaces. Recall that a Banach space is {\em weakly compactly generated} (briefly: {\em WCG}) if it contains a linearly dense weakly compact set.
\index{WCG space}\index{Banach space!-- WCG}
\index{Banach space!-- weakly compactly generated}

\begin{prop}\label{iowetjruqw}
Let $X$ be a weakly compactly generated Banach space and let $\theta$ be a big enough regular cardinal. Further, let $M\rloe \pair{H(\theta)}\in$ be such that $X\in M$. Then there exists a norm one projection $\map {P_M}X{X_M}$ such that $\ker(P_M) = \leftorth{X^*\cap M}$.
\end{prop}

\begin{pf}
Let $K$ be a linearly dense weakly compact subset of $X$. By Lemma \ref{sgagaqqwfwr}, it suffices to check that $X_M \cup \leftorth{X^*\cap M}$ is linearly dense in $X$.

Suppose $\phi\in X^*\setminus\sn0$ is such that $(X\cap M)\subs\ker(\phi)$ and $\leftorth{X^*\cap M}\subs \ker(\phi)$. The latter inclusion implies that $\phi \in \clstar(X^*\cap M)$, because $X^*\cap M$ is $\Qyu$-linear. Fix $p\in K$ such that $\phi(p)\ne 0$. Let $U_0, U_1\subs \Err$ be disjoint open rational intervals such that $0\in U_0$ and $\phi(p)\in U_1$.
Let $K_0$ be the weak closure of $K\cap M$. Note that $\phi\rest K_0=0$, because $\phi$ is weakly continuous. Using the fact that $\phi \in \clstar(X^*\cap M)$, for each $x\in K_0$ choose $\psi_x\in X^*\cap M$ such that $\psi_x(x)\in U_0$ and $\psi_x(p)\in U_1$. By compactness, there are $x_0,x_1,\dots,x_{n-1}\in K_0$ such that
\begin{equation}
K_0\subs \bigcup_{i<n}\psi_{x_i}^{-1}[U_0].
\tag{$*$}
\end{equation}
Note that $U_0,U_1\in M$. Let $\Psi=\setof{\psi_{x_i}}{i<n}$. By Proposition \ref{enbgojhe}(b), also $\Psi\in M$. Further, $p\in K$ witnesses the validity of
$$(\exists\;x\in K)(\forall\;\psi\in\Psi)\;\;\psi(x)\in U_1.$$
All parameters in the above formula are in $M$, therefore by elementarity there exists $x\in K\cap M\subs K_0$ such that $\psi_{x_i}(x)\in U_1$ for each $i<n$. This contradicts ($*$).
\end{pf}

\index{norming space/set!-- generating projections}
Fix a Banach space $X$ and a norming set $D\subs X^*$. We shall say that $D$ {\em ge\-ne\-rates projections in} $X$ if for a big enough regular cardinal $\theta$, for every countable elementary substructure $M\rloe\pair{H(\theta)}\in$ with $D\in M$ it holds that
$$X=X_M\oplus\leftorth{D\cap M}.$$
We shall say that the projection $\map PXX$ such that $\im P=X_M$ and $\ker P=\leftorth{D\cap M}$ is {\em induced by} the triple $\triple XDM$.
\index{projection!-- induced by $\triple XDM$}
Proposition \ref{iowetjruqw} says that ${X^*}$ generates projections in $X$ whenever $X$ is WCG.
The class of Banach spaces $X$ with the property that ${X^*}$ generates projections in $X$ is well known: these are precisely weakly Lindel\"of determined spaces. Below we prove the easier implication.

\begin{prop}\label{poopipriuq}
Let $\pair XD$ be a Plichko pair. Then $D$ generates projections in $X$.
\end{prop}

A similar statement was proved by Koszmider \cite[Lemma 4.1]{Kosz-wcg}, using the existence of a countably $1$-norming Markushevich basis.

\begin{pf}
Fix $M\rloe \pair{H(\theta)}\in$ such that $D\in M$.
Suppose $\phi\in X^*$ is such that $X\cap M\subs \ker \phi$ and $\leftorth{D\cap M}\subs \ker\phi$. The latter fact means that $\phi\in\clstar(D\cap M)$. Let $G\subs X$ be a linearly dense set such that $\suppt(y,G)$ is countable for each $y\in D$. By elementarity, we may assume that $G\in M$. 
Suppose $\phi\ne0$ and fix $u\in G$ such that $\phi(u)\ne0$. Since $\phi$ is in the \weakstar\ closure of $D\cap M$, we may find $\psi\in D\cap M$ such that $\psi(u)\ne0$. Thus $u\in\suppt(\psi,G)$. On the other hand, $\suppt(\psi,G)\in M$, because $\psi,G\in M$. By Proposition \ref{enbgojhe}(c), $\suppt(\psi,G)\subs M$. In particular $u\in X\cap M$ and hence $\phi(u)=0$, a contradiction.
\end{pf}

The above result says in particular that ${X^*}$ generates projections in $X$ whenever $X$ is WLD. The converse implication will be proved in Section \ref{popop}.

\subsection{Projectional generators}

Projectional generators were introduced by Orihuela and Valdivia \cite{OriVal90} as a tool for showing that certain Banach spaces have a PRI. In fact, first projectional generators were implicitely constructed by Lindenstrauss \cite{Lin65, Lin66}. We refer to Chapter 6 of Fabi\'an's book \cite{Fab97} for more information.
Let us recall the definition.

\index{projectional generator}
Let $X$ be a Banach space.
A pair $\pair D\Phi$ is a {\em projectional generator} in $X$ if
\begin{enumerate}
	\item[(1)] $D\subs X^*$ is a norming $\Qyu$-linear subspace,
	\item[(2)] $\map\Phi D{\dpower X{\loe\aleph_0}}$,
	\item[(3)] $(\bigcup\img\Phi B)^\perp\cap \clstar(B)=0$, whenever $B\subs D$ is $\Qyu$-linear.
\end{enumerate}
Below we show how projectional generators together with elementary submodels induce projections.

\begin{prop}
Let $X$ be a Banach space which has a projectional generator with domain $D$. Then $D$ generates projections in $X$.
\end{prop}

\begin{pf}
Fix $M\rloe\pair{H(\theta)}\in$ with $D\in M$, where $\theta$ is big enough and, using elementarity, assume $\Phi\in M$ is a projectional generator with domain $D$.
It suffices to check that the only $\psi\in X^*$ which vanishes on $X_M\cup\leftorth{D\cap M}$ is the zero functional. 

Let $B:=D\cap M$. By elementarity, $B$ is $\Qyu$-linear. By the definition of a projectional generator, $(\bigcup\img\Phi B)^\perp\cap \clstar(B)=0$. Thus, if $\psi\in X^*$ is such that $X_M\cup\leftorth{D\cap M}\subs \ker\psi$ then $\psi\in (\bigcup\img\Phi B)^\perp$. This is because $\Phi(b)\subs M$ whenever $b\in B$ (by Proposition \ref{enbgojhe}(c)). It follows that $\psi=0$.
\end{pf}

\subsection{Bandlow's Property $\Omega$}

A result of Bandlow \cite[Thm. 5.6]{Ban94} says that $K$ is Corson compact iff $\ceep K$ has Property $\Omega$, which says that for every sufficiently closed countable $M\rloe\pair{H(\theta)}\in$, for every $f\in \ceep K$ there exists $g\in\cl(\ceep K\cap M)$ such that $f\rest (K\cap M) = g\rest (K\cap M)$ (in other words: $f-g\in\leftorth{K\cap M}$).
According to Bandlow's definition, $\cl$ means here the pointwise closure, however by Lemma \ref{ergnosfoipq3r} we can replace it by the norm closure.

\index{Property $\Omega$}
A natural generalization of the above condition is {\em Property $\Omega$} for a pair $\pair XD$, where $D$ is a norming set in the dual of $X$: given a suffciently closed countable elementary substructure $M$ of $H(\theta)$, for every $p\in X$ there is $q\in \cl(X\cap M)$ such that $p-q\in\leftorth{D\cap M}$.

\begin{prop}\label{weriojqwpriq}
Let $X$ be a Banach space and let $D\subs X^*$ be a norming set. The following statements are equivalent.
\begin{enumerate}
	\item[(a)] $\pair XD$ has Property $\Omega$.
	\item[(b)] $D$ generates projections in $X$.
\end{enumerate}
\end{prop}

\begin{pf}
(a)$\implies$(b) Fix a suitable $M\rloe\pair{H(\theta)}\in$ and fix $\phi\in X^*$ such that $(X\cap M)\cup\leftorth{D\cap M}\subs \ker \phi$. Fix $p\in X$. Applying Property $\Omega$, find $q\in \cl(X\cap M)$ such that $p-q\in\leftorth{D\cap M}$. Then $\phi(q)=0$ and $\phi(p-q)=0$, thus also $\phi(p)=0$. This shows that $\phi=0$. Since $\phi$ was arbitrary, we get $X=\cl(X\cap M)\oplus\leftorth{D\cap M}$.

(b)$\implies$(a) Fix a suitable $M\rloe\pair{H(\theta)}\in$ and fix $p\in X$. By (b), $X=\cl(X\cap M)\oplus\leftorth{D\cap M}$, therefore there is $q\in \cl(X\cap M)$ such that $p-q\in\leftorth{D\cap M}$. This is exactly Property $\Omega$.
\end{pf}

Thus, in our language, Bandlow's result says that $K$ is Corson compact if and only if $K$ generates projections in $\cee K$, where $K$ is naturally identified with a suitable subset of $\cee K^*$. Recall that consistently there exists a Corson compact $K$ for which $\cee K$ is not WLD, see \cite{ArgMerNeg}.

\section{Projectional skeletons}

In this section we define the crucial notion of this work.
Recall that a partially ordered set $\Gam$ is {\em directed} if for every $s_0,s_1\in \Gam$ there is $t\in\Gam$ such that $s_0\loe t$ and $s_1\loe t$. $\Gam$ is {\em $\sig$-complete} if every sequence $s_0<s_1<\dots$ has the least upper bound in $\Gam$. A subset $A$ of $\Gam$ is {\em closed} if $\sup_{\ntr} s_n\in A$ whenever $\setof{s_n}{\ntr}\subs A$ is such that $s_0<s_1<\dots$. A set $A\subs \Gam$ is {\em cofinal} if for every $s\in \Gam$ there exists $t\in A$ with $s\loe t$.

\subsection{Definition and basic properties}

\index{projectional skeleton}
Let $X$ be a Banach space. A {\em projectional skeleton} in $X$ is a family $\sett{P_s}{s\in\Gam}$ of bounded projections of $X$ indexed by a directed partially ordered set $\Gam$, satisfying the following conditions
\begin{enumerate}
	\item[(1)] $X=\bigcup_{s\in\Gam}P_s X$ and each $P_s X$ is separable.
	\item[(2)] $s\loe t\implies P_s=P_s\cmp P_t = P_t \cmp P_s$.
	\item[(3)] If $s_0<s_1<s_2<\dots$ then $t=\sup_{\ntr}s_n$ exists in $\Gam$ and $P_t X=\cl(\bigcup_{\ntr}P_{s_n} X)$.
\end{enumerate}
Condition (3) says in particular that the poset $\Gam$ is $\sig$-complete.
We have not assumed so far that the projections $P_s$ are uniformly bounded. Note that for every closed cofinal set $\Gam'\subs \Gam$ the restriction $\sett{P_s}{s\in\Gam'}$ is again a projectional skeleton in $X$.
The notion of a projectional skeleton makes sense for non-separable Banach spaces only: in a separable Banach space $X$ the family $\sn{\id X}$ is a projectional skeleton.

The next observation was obtained jointly with Ond\v rej Kalenda.

\begin{prop}
Let $\sett{P_s}{s\in\Gam}$ be a projectional skeleton in a Banach space $X$. Then there exists a closed cofinal set $\Gam'\subs \Gam$ such that
$$\sup_{s\in\Gam'}\norm{P_s} < +\infty.$$
\end{prop}

\begin{pf} For each $n\goe1$ define $G_n=\setof{s\in\Gam}{\norm{P_s}\loe n}$. We claim that one of these sets is cofinal in $\Gam$. Suppose otherwise and for each $\ntr$ choose $t_n$ such that $\norm{P_s}>n$ whenever $t_n\loe s$. Using the directedness of $\Gam$, construct a sequence $s_1<s_2<\dots$ such that $t_n\loe s_n$ for $\ntr$. Let $s_\infty=\sup_{\ntr}s_n$. Then $\norm{P_{s_\infty}}=+\infty$, a contradiction.

Fix $k\goe1$ such that $\Gam':=G_k$ is cofinal in $\Gam$. We claim that $\Gam'$ is also closed. For fix $s_0<s_1<\dots$ in $\Gam'$ and let $t=\sup_{\ntr}s_n$. We need to show that $\norm{P_t}\loe k$.

Suppose this is not true and fix $x\in X$ with $\norm x=1$ and $\norm{P_t x}=r>k$. Let $\eps=(r-k)/2$. Using the second part of (3), find $m\in\nat$ and $y\in P_{s_m}X$ such that $\norm{P_t x-y}<\eps/k$. Note that
$P_{s_m}=P_{s_m}\cmp P_t$. Using the fact that $\norm{P_{s_m}}\loe k$, we get $\norm{P_{s_m}x}\loe k$ and
$$\norm{y - P_{s_m}x} =\norm{P_{s_m}(y - P_t x)} \loe k\norm{y - P_t x}<\eps.$$
Thus
$$
\norm{P_tx} \loe \norm{P_tx - y} +\norm{y - P_{s_m}x} + \norm{P_{s_m}x} < \eps/k+\eps+k \loe 2\eps+k = r = \norm{P_t x}.
$$
This contradiction completes the proof.
\end{pf}

By the above proposition, we shall always assume that a projectional skeleton $\sett{P_s}{s\in\Gam}$ satisfies the condition
\begin{enumerate}
	\item[(4)] $\sup_{s\in\Gam}\norm{P_s} < +\infty$.
\end{enumerate}
We shall say that $\sett{P_s}{s\in\Gam}$ is an {\em $r$-projectional skeleton} if it is a projectional skeleton such that $\norm{P_s}\loe r$ for every $s\in\Gam$.
The remaining part of this section will be devoted to proving basic properties of projectional skeletons.

\begin{lm}\label{eiurqrad}
Let $\sett{P_s}{s\in\Gam}$ be a projectional skeleton in $X$ and let $s_0<s_1<\dots$ be such that $t=\sup_{\ntr}s_n$ in $\Gam$. Then
$$P_tx=\lim_{n\to\infty}P_{s_n}x$$
for every $x\in X$.
\end{lm}

\begin{pf}
Let $r=\sup_{s\in\Gam}\norm{P_s}$ and fix $x\in X$, $\eps>0$. By the second part of (3), find $y\in \bigcup_{\ntr}P_{s_n}X$ such that $\norm{P_tx - y} < \eps/(1+r)$. Choose $k$ such that $y\in P_{s_{k}}X$. Note that $P_ty=y$ and $P_{s_n}y=y$ for $n\goe k$. Thus, given $n\goe k$, we have
\begin{align*}
\norm{P_tx - P_{s_n}x}&\loe \norm{P_tx - y} + \norm{y - P_{s_n}x} < \eps/(1+r) + \norm{P_{s_n}(y - P_tx)}\\
&\loe \eps/(1+r) + r\norm{y - P_tx} < \eps/(1+r) + r\eps/(1+r) = \eps.
\end{align*}
This shows that $\lim_{n\to\infty}\norm{P_tx-P_{s_n}x}=0$.
\end{pf}

\begin{lm}\label{iuwbfuiq}
Let $\sett{P_s}{s\in\Gam}$ be a projectional skeleton in $X$ and let $T\subs\Gam$ be a directed subset of $\Gam$. Then the formula
$$P_Tx = \lim_{s\in T}P_sx$$
well defines a bounded projection of $X$ onto $\cl(\bigcup_{s\in T}P_sX)$.
\end{lm}

\begin{pf}
We must show that $\sett{P_sx}{s\in T}$ is a Cauchy net for every $x\in X$.
Fix a big enough regular cardinal $\theta$ so that all relevant objects are contained in $H(\theta)$. Fix a countable elementary substructure $M$ of $H(\theta)$ such that $\sett{P_s}{s\in\Gam}\in M$ and $T\in M$. Let $S=T\cap M$. By elementarity, $S$ is a directed set. Suppose that there exists $x\in X$ such that $\sett{P_sx}{s\in T}$ is not a Cauchy net. By elementarity, there exist $x\in M\cap X$ and $\eps>0$ in $M$ such that
$$M\models \forall\;t\in T\;\exists\;s,s'\in T\;(s\goe t\land s'\goe t\land \norm{P_sx - P_{s'}x}\goe\eps).$$
This means in particular that, given $t\in S$, there exist $s,s'\in S$ such that $s,s'\goe t$ and $\norm{P_sx - P_{s'}x}\goe\eps$. In order to get a contradiction, it suffices to show that $\sett{P_sx}{s\in S}$ is a Cauchy net.

Since $S$ is countable and directed, it has an increasing cofinal sequence $s_0<s_1<\dots$. Let $u=\sup_{\ntr}s_n$. Then $u=\sup S$. By Lemma \ref{eiurqrad}, $P_tx=\lim_{n\to\infty}P_{s_n}x$ for every $x\in X$. Fix $\eps>0$ and find $k$ such that $\norm{P_{s_n}x - P_{s_{k}}x}<\eps/(2r)$ for every $n\goe k$, where $r\goe1$ is such that $\norm{P_s}\loe r$ for every $s\in T$.
Fix $s,t\in S$ such that $s,t\goe s_{k}$. Choose $\ell>k$ so that $s,t\loe s_{\ell}$. We have
$$\norm{P_sx - P_{s_k}x} = \norm{P_sP_{s_\ell}x - P_sP_{s_k}x}\loe\norm{P_s}\cdot\norm{P_{s_\ell}x - P_{s_k}x}<r\eps/(2r) = \eps/2.$$
Similarly, $\norm{P_tx - P_{s_k}x}<\eps/2$. Thus, finally $\norm{P_sx-P_tx}<\eps$ for every $s,t\in S$ with $s,t\goe s_k$, i.e. $\sett{P_sx}{s\in S}$ is a Cauchy net for every $x\in X$. 

This shows that $P_T$ is well defined. It is clear that $P_T$ is a bounded projection and that $P_Tx=x$ iff $x\in \cl(\bigcup_{s\in T}P_sX)$.
\end{pf}

\subsection{Projectional resolutions of the identity}

\begin{tw}\label{qwrorqw}
Every Banach space with $1$-projectional skeleton has a projectional resolution of the identity.
\end{tw}

\begin{pf}
Let $\kappa=\dens X$ and let $\s = \sett{P_s}{s\in\Gam}$ be a projectional skeleton in $X$ such that $\norm{P_s}=1$ for every $s\in\Gam$. Fix a continuous chain $\sett{T_{\al}}{\al<\kappa}$ of up-directed subsets of $\Gam$ satisfying $|T_\al|\loe\al+\aleph_0$ for each $\al$ and such that
$$E=\bigcup\setof{P_s X}{s\in T_\al,\; \al<\kappa}$$
is dense in $X$. Continuity of the chain means that $T_\delta=\bigcup_{\xi<\delta}T_\xi$ whenever $\delta$ is a limit ordinal.
Let $X_\al=\cl(\bigcup_{s\in T_\al}P_sX)$.
Then $\sett{X_\al}{\al<\kappa}$ is a chain of closed subspaces of $X$, the density of $X_\al$ does not exceed $|\al|+\aleph_0$ and $X_\delta = \cl(\bigcup_{\xi<\delta}X_\xi)$ for every limit ordinal $\delta$.
By Lemma \ref{iuwbfuiq}, formula $P_\al x = \lim_{s\in T_\al}P_sx$ defines a projection of $X$ onto $X_\al$.
Clearly, $\norm{P_\al}=1$, because $\norm{P_s}=1$ for $s\in \Gam$. It is easy to check, again using Lemma \ref{iuwbfuiq}, that $P_\al = P_\al\cmp P_\beta = P_\beta\cmp P_\al$ for $\al<\beta$.
Thus, $\sett{P_\al}{\al<\kappa}$ is a PRI on $X$.
\end{pf}

\begin{wn}
Given a Banach space $X$ of density $\aleph_1$, the following properties are equivalent.
\begin{enumerate}
	\item[(a)] $X$ has a $1$-projectional skeleton.
	\item[(b)] $X$ has a projectional resolution of the identity.
\end{enumerate}
\end{wn}

\begin{pf}
Implication (a)$\implies$(b) follows from the above theorem. In case of density $\aleph_1$, every PRI is a $1$-projectional skeleton.
\end{pf}

\subsection{Norming space induced by a projectional skeleton}

We shall now look at the dual of a space with a projectional skeleton.
Let $\s = \sett{P_s}{s\in\Gam}$ be a projectional skeleton in a Banach space $X$.
The set
$$D = \bigcup_{s\in\Gam}P_s^*X^*$$
is clearly a linear subspace of $X$. Notice that $P_s^*X^*\cap\ubal{X^*}$ endowed with the \weakstar\ topology is second countable, because $P_s^*X^*$ is linearly homeomorphic to the dual of $P_sX$.
Let $r=\sup_{s\in\Gam}\norm{P_s}$. Given $x\in\usphere X$, there is $s\in\Gam$ such that $x=P_sx$; choose $\phi\in X^*$ such that $\phi(x)=1=\norm\phi$. Then $(P_s^*\phi)x = \phi(P_s x) = \phi(x)=1$ and $\norm{P_s^*\phi}\loe r$. This shows that the space $D$ is $r$-norming. We shall say that $D$ is the dual norming subspace {\em induced by} $\s$ and we shall denote it by $D(\s)$.
\index{$D(\argum)$}

\index{projectional skeleton!-- equivalence of}
Let $\s=\sett{P_s}{s\in\Gam}$ and $\ft=\sett{Q_t}{t\in\Pi}$ be projectional skeletons in the same Banach space $X$. We say that $\s$ and $\ft$ are {\em equivalent} if they induce the same norming subspace, i.e. $\bigcup_{s\in\Gam}\im P_s^* = \bigcup_{t\in\Pi}\im Q_t^*$.
It turns out that, with help of elementary submodels, a projectional skeleton can be recovered (up to equivalence) from the norming space.

\begin{lm}\label{weiopfjpqw}
Let $\s=\sett{P_s}{s\in\Gam}$ be a projectional skeleton in a Banach space $X$ and let $D\subs D(\s)$ be norming for $X$. Further, let $\theta$ be a big enough regular cardinal and let $M\rloe \pair{H(\theta)}\in$ be countable and such that $\s\in M$. Then the projection induced by $\triple XDM$ equals $P_t$, where $t=\sup(\Gam\cap M)$.
\end{lm}

\begin{pf} 
Observe that
$$\bigcup_{s\in \Gam\cap M}\im P_s\subs \cl(X\cap M).$$
This is because, given $s\in \Gam\cap M$, by elementarity there exists a countable set $A\in M$ which is dense in $\im P_s$. By Proposition \ref{enbgojhe}(c), $A\subs X\cap M$, therefore $\im P_s\subs \cl(X\cap M)$. It follows that $\im P_t = \cl(\bigcup_{s\in \Gam\cap M}\im P_s)\subs \cl(X\cap M)$. On the other hand, given $x\in X\cap M$, by elementarity there is $s\in \Gam\cap M$ such that $x\in\im P_s$; thus $x\in \im P_t$. Hence $\im P_t = X_M$.
Notice that, again by elementarity, $P_t^*\phi=\phi$ whenever $\phi\in D\cap M$. Thus $\ker P_t\subs \leftorth{D\cap M}$, because if $P_tx=0$ and $\phi\in D\cap M$ then $\phi(x)=(P_t^*\phi)x=\phi(P_tx)=0$. It follows that $\ker P_t = \leftorth{D\cap M}$, because $X=X_M\oplus\ker P_t$ and $X_M\cap\leftorth{D\cap M}=\sn0$ (Lemma \ref{sgagaqqwfwr}(a)).
\end{pf}

\begin{tw}\label{wejirpqwrt}
Let $X$ be a Banach space and let $D\subs X^*$ be an $r$-norming set ($r\goe1$). The following properties are equivalent.
\begin{enumerate}
	\item[(a)] $X$ has an $r$-projectional skeleton $\s$ such that $D\subs D(\s)$.
	\item[(b)] $D$ generates projections in $X$.
\end{enumerate}
\end{tw}

\begin{pf}
Implication (a)$\implies$(b) follows from Lemma \ref{weiopfjpqw}.

(b)$\implies$(a)
Fix a big enough regular cardinal $\theta$ and let $\Gam$ be the family of all countable elementary substructures $M$ of $\pair{H(\theta)}\in$ such that $D\in M$. Endow $\Gam$ with inclusion. Clearly, $\Gam$ is a $\sig$-directed poset. Fix 
$M\in \Gam$ and let $P_M$ be the projection onto $\cl(X\cap M)$ with $\ker(P_M)=\leftorth{D\cap M}$. By definition, $P_M$ is defined on the entire space.  Further, $\norm{P_M}\loe r$, by Lemma \ref{sgagaqqwfwr}.
Given a sequence $M_0\subs M_1\subs\dots$ in $\Gam$, the union $M=\bigcup_{\ntr}M_n$ is an elementary substructure of $H(\theta)$ such that $\cl(X\cap M)$ is the closure of $\bigcup_{\ntr}\cl(X\cap M_n)$. It follows that $\s=\sett{P_M}{M\in\Gam}$ is a projectional skeleton in $X$. It is clear that $D\subs D(\s)$, because $D\cap M\subs\im P_M^*$.
\end{pf}

\begin{wn}\label{werour}
Let $X$ be a Banach space with a projectional skeleton. Then there exists a renorming of $X$ under which $X$ has an equivalent $1$-projectional skeleton.
\end{wn}

\begin{pf}
Let $D$ be the dual norming subspace induced by a fixed projectional skeleton in $X$. By Lemma \ref{weiopfjpqw}, $D$ generates projections in $X$. Consider a renorming of $X$ after which $D$ becomes $1$-norming (see Proposition \ref{nofjsopfj}). By Theorem \ref{wejirpqwrt}, $X$ has a $1$-projectional skeleton.
\end{pf}

As an application of Theorem \ref{wejirpqwrt}, we prove that the class of Banach spaces with a $1$-projectional skeleton is stable under arbitrary $c_0$- and $\ell_p$-sums.

\begin{tw}\label{eoghjnoasf}
Let $\sett{X_\al}{\al<\kappa}$ be a collection of Banach spaces and let $X=\bigoplus_{\al<\kappa}X_\al$ be endowed either with the $c_0$-norm or with $\ell_p$-norm ($1\loe p <\infty$). Further, assume that for each $\al<\kappa$, $D_\al\subs X_\al$ is $1$-norming and generates projections in $X_\al$. Then the set
$$D=\setof{\phi\in X^*}{(\forall\;\al)\;\phi\rest X_\al\in D_\al\;\text{ and }\;|\setof{\al}{\phi\rest X_\al\ne0}|\loe\aleph_0}$$
is $1$-norming and generates projections in $X$.
\end{tw} 

\begin{pf}
The fact that $D$ is $1$-norming follows from the properties of the $c_0$-sum and the $\ell_p$-sum.
Define $\suppt(\phi) = \setof{\al}{\phi\rest X_\al\ne0}$. Fix a suitable $M\rloe\pair{H(\theta)}\in$ so that $D\in M$. Let $S=\kappa\cap M$. Note that $\suppt(\phi)\subs S$ whenever $\phi\in D\cap M$. Suppose $X\ne X_M\oplus\leftorth{D\cap M}$ and fix $\psi\in X^*$ satisfying $X\cap M\subs\ker \psi$ and $\leftorth{D\cap M}\subs \ker \psi$. Then $\psi$ is in the \weakstar\ closure of the linear hull of $D\cap M$, therefore $\suppt(\psi)\subs S$. Assuming $\psi\ne0$, there is $\al\in S$ such that $\psi_\al:=\psi\rest X_\al\ne0$. Note that $X_\al\cap M\subs \ker \psi_\al$. If $x\in \leftorth{D_\al\cap M}$ and $\phi\in D\cap M$ then $\phi\rest X_\al\in D_\al\cap M$ so $\phi(x)=0$. It follows that $\leftorth{D\cap M}\cap X_\al = \leftorth{D_\al\cap M}$. Thus $\leftorth{D_\al\cap M}\subs \ker \psi$. On the other hand, $X_\al = \cl(X_\al\cap M)\oplus \leftorth{D_\al\cap M}$, because $D_\al\in M$. This is a contradiction.
\end{pf}

We finish this section by exhibiting a topological property of norming spaces induced by a projectional skeleton.

\begin{tw}\label{rwqroiufd}
Let $X$ be a Banach space with a projectional skeleton $\s=\sett{P_s}{s\in\Gam}$ and let $D\subs X^*$ be the norming space induced by $\s$, endowed with the \weakstar\ topology. Then:
\begin{enumerate}
	\item[(a)] The closure in $X^*$ of every countable bounded subset of $D$ is metrizable and contained in $D$.
	\item[(b)] $D$ is countably tight.
\end{enumerate}
\end{tw}

\begin{pf}
Part (a) is trivial: every countable subset of $D$ is contained in $Y_s=P^*_sX^*$ for some $s\in\Gam$ and every bounded subset of $Y_s$ with the \weakstar\ topology is second countable.

(b) Let $A\subs D$ and $p\in\clstar(A)\cap D$ be given. Replacing $A$ by $A-p$, we may assume that $p=0$. Fix a big enough regular cardinal $\theta$ and a countable elementary substructure $M$ of $\pair{H(\theta)}\in$ such that $X,\s,A\in M$. We claim that $0\in\clstar(A\cap M)$.

Let $t=\sup(\Gam\cap M)$ and let $Y_t=P_t^*X^*$. Fix a $\weakstar$ neighborhood $U$ of $p$. We may assume that $U=\bigcap_{i<n}U(x_i,\eps)$, where $U(x,\eps) := \setof{y\in X^*}{|y(x)|<\eps}$, $x_0,\dots,x_{n-1}\in X$ and $\eps>0$ is rational.
By Lemma \ref{weiopfjpqw}, $P_tX=\cl(X\cap M)$. By Banach's Open Mapping Principle, $\inv{P_t}{X\cap M}$ is dense in $X$. Hence, without loss of generality, we may assume that $P_tx_i\in M$ for each $i<n$.

Thus $W:=\bigcap_{i<n}U(P_tx_i,\eps)$ is a \weakstar\ neighborhood of $0$ and $W\in M$, because $P_tx_i\in M$ and $\eps\in M$. By elementarity, there is $a\in A\cap M\cap W$. It follows that $a\in Y_t$, i.e. $P_t^*a=a$. Given $i<n$, we have $a(x_i) = (P_t^*a)x_i = a(P_tx_i)$. Thus $a\in U(x_i,\eps)$. Finally, $a\in A\cap M\cap U$.
\end{pf}

\begin{wn}
Let $X$ be a Banach space and let $D,E\subs X^*$ be norming spaces induced by projectional skeletons. If $D\cap E$ is total then $D=E$.
\end{wn}

\begin{pf}
Note that $D\cap E$ is a linear space. Assuming it is total, it must be \weakstar\ dense. Thus, given $p\in D$, we have that $p\in\clstar(D\cap E)$ so $p\in \clstar(A)$ for some countable $A\subs D\cap E$ (Theorem \ref{rwqroiufd}(b)). By Theorem \ref{rwqroiufd}(a), $\clstar(A)\subs D\cap E$. This shows that $D\subs E$. By symmetry, $D=E$. 
\end{pf}

\begin{wn}
Let $X$ be a Banach space and let $S\subs X^*$ generate projections in $X$. Further, let $D\sups S$ be the smallest linear subspace of $X^*$ such that $\clstar(A)\subs D$ for every countable set $A\subs D$. Then $D$ is induced by a projectional skeleton in $X$.
\end{wn}

\begin{pf}
By Theorem \ref{wejirpqwrt}, there exists a projectional skeleton $\s$ in $X$ such that $S\subs D(\s)$. By Theorem \ref{rwqroiufd}(a), $D\subs D(\s)$. On the other hand, $D$ is \weakstar\ dense, because it is a norming (and hence total) linear space. Hence, given $p\in D(\s)$, we have that $p\in \clstar(D)$ so, by Theorem \ref{rwqroiufd}(b), $p\in \clstar(A)$ for some countable $A\subs D$. Hence $p\in D$. This shows that $D=D(\s)$.
\end{pf}

The results of this section show that Banach spaces with a projectional skeleton have very similar properties to Plichko spaces. In fact, we know only one basic example distinguishing those two classes: the space $\cee{\omega_2+1}$ which, by the result of Kalenda \cite{Kal02}, does not have any countably norming Markushevich basis.

\section{Plichko spaces and projectional skeletons}\label{popop}

We prove a preservation theorem for projectional sequences of Plichko spaces. As an application, we show that every Banach space with a commutative projectional skeleton is Plichko.
\index{projectional skeleton!-- commutative}
A projectional skeleton $\sett{P_s}{s\in\Gam}$ is {\em commutative} if $P_s\cmp P_t=P_t\cmp P_s$ holds for every $s,t\in\Gam$.

\begin{prop}\label{ofhwoehfr}
Let $X$ be a Plichko space with a $\Sigma$-space $D\subs X^*$. Then there exists a commutative projectional skeleton $\sett{P_s}{s\in\Gam}$ on $X$ such that $D=D(\s)$.
\end{prop}

\begin{pf}
Fix a big enough regular cardinal $\theta$ and let $\Gam$ be the family of all countable elementary substructures $M$ of $\pair{H(\theta)}\in$ such that $\pair XD\in M$. By Proposition \ref{poopipriuq} and the proof of Theorem \ref{wejirpqwrt}, we know that $\s=\sett{P_M}{M\in\Gam}$ is a projectional skeleton in $X$ such that $D\subs D(\s)$, where $P_M$ is induced by $\triple XDM$, i.e. $\im P_M=X_M$ and $\ker P_M=\leftorth{D\cap M}$. Theorem \ref{rwqroiufd}(b) says that $D(\s)$ is \weakstar\ countably tight. On the other hand, $D$ is \weakstar\ countably closed, i.e. $\clstar A\subs D$ whenever $A\subs D$ is countable. Thus $D=D(\s)$.

It remains to show that $\s$ is commutative. Given $M\in\Gam$, define $r_M = P_M\rest G$, where $G\subs X$ is a linearly dense set witnessing that $D$ is a $\Sigma$-space. We may assume that $G\in M$.
We claim that
\begin{equation}
r_M(x)=
\begin{cases}
x,\quad&\text{ if }x\in G\cap M,\\
0,&\text{ if }x\in G\setminus M.
\end{cases}
\tag{*}
\end{equation}
Indeed, $G\cap M\subs \cl(X\cap M)=\im P_M$. If $x\in G\setminus M$ then for $y\in D\cap M$ we have that $x\notin\suppt(y,G)$, because $\suppt(y,G)\subs M$. Hence $y(x)=0$ for $y\in D\cap M$ and therefore $x\in\leftorth{D\cap M}=\ker P_M$. 

Using (*), we see that $r_M\cmp r_N=r_{M\cap N}$ for $M,N\in\Gam$.
Since $G$ is linearly dense in $X$, this shows that $P_M\cmp P_N = P_{M\cap N} = P_N\cmp P_M$ for every $M,N\in\Gam$. This completes the proof.
\end{pf}

\subsection{Preservation theorem}

\begin{lm}\label{fasejqwor}
Let $\pair XD$ be a Plichko pair and let $\map PXX$ be a bounded projection such that $P^*D\subs D$. Then $\pair{\ker P}{D\cap \ker P^*}$ is a Plichko pair.
\end{lm}

\begin{pf} Let $A$ be a linearly dense subset of $X$ such that $|\suppt(y,A)|\loe\aleph_0$ for every $y\in D$.
Let $B=\setof{a-Pa}{a\in A}$. Then $B$ is linearly dense in $\ker P = \im(\id X - P)$. Let $E=D\cap \ker P^*$. Fix $y\in D$ with $\norm y=1$. Let $z=y-P^*y$. Then $z\in D$, because $P^*$ preserves $D$. Further, $z\in \ker P^*$ and $\norm z\loe 1+\norm P$. Given $x\in \ker P$, we have $|z(x)|=|y(x)-(P^*y)x|=|y(x)-y(Px)|=|y(x)|$.
This shows that $E$ is norming for $\ker P$. Finally, given $y\in E$, we have $y(a-Pa)=y(a)-y(Pa)=y(a)-(P^*y)a=y(a)$, so $\suppt(y,B)=\suppt(y,A)$.
Thus, $B$ witnesses that $\pair{\ker P}{E}$ is a Plichko pair.
\end{pf}

\begin{tw}\label{popiooier}
Let $\sett{P_\al}{\al<\kappa}$ be a projectional sequence in a Banach space $X$ and let $D\subs X^*$ be a norming space such that
$$D=\bigcup_{\al<\kappa}P_\al^* D$$
and $\pair{P_\al X}{P_\al^* D}$ is a Plichko pair for each $\al<\kappa$.
Then $\pair XD$ is a Plichko pair.
\end{tw}

Note that we do not assume that the above projections are uniformly bounded.

\begin{pf}
We construct inductively a family of sets $\sett{A_\al}{\al<\kappa}$ such that
\begin{enumerate}
	\item[(i)] $A_\al$ is a linearly dense subset of $P_\al X$;
	\item[(ii)] $\al<\beta\implies A_\al\subs A_\beta$;
	\item[(iii)] $\suppt(y,A_\al)$ is countable for every $y\in D$;
	\item[(iv)] $P_\al a=0$ whenever $a\in A_\beta\setminus A_\al$.
\end{enumerate}
We start with a linearly dense set $A_0\subs P_0X$ witnessing that $\pair{P_0X}{P_0^*D}$ is a Plichko pair. We must check (iii). Observe that $a=P_0a$ and consequently $y(a)=(P_0^*y)a$ for every $a\in A_0$ and $y\in X^*$. Hence $\suppt(y,A_0)=\suppt(P_0^*y,A_0)$ is countable and (iii) holds.
Now, fix an ordinal $\delta>0$ and assume $\sett{A_\al}{\al<\delta}$ has already been defined.

Suppose first that $\delta=\beta+1$. We use Lemma \ref{fasejqwor}. 
Let $Y=P_\delta X$ and let $D'=\setof{y\rest Y}{y\in D}\subs Y^*$. Then $\pair Y{D'}$ is a Plichko pair and $P_\beta\rest Y$ is a projection whose dual preserves $D'$. By Lemma \ref{fasejqwor}, there is a linearly dense subset $B$ of $\ker P_\beta\cap Y$ witnessing that $\pair {\ker P_\beta\cap Y}E$ is a Plichko pair, where $E=D'\cap \ker [(P_\beta\rest Y)^*]$. It follows that $\suppt(y,B)$ is countable whenever $y\in D\cap\ker P_\beta^*$.
Define $A_\delta = A_\beta\cup B$. Conditions (i), (ii) and (iv) are obviously satisfied. It remains to check (iii). Fix $y\in D$. By the induction hypothesis, $\suppt(y,A_\beta)$ is countable. Let $z=y-P_\beta^*y$. Then $z\in D\cap\ker P_\beta$, so $\suppt(z,B)$ is countable. Finally, given $b\in B$ we have
$$z(b) = y(b) - (P_\beta^*y)b = y(b) - y(P_\beta b) = y(b),$$
therefore $\suppt(y,B)=\suppt(z,B)$ is countable. This shows (iii), because $\suppt(y,A_\delta)$ $=\suppt(y,A_\beta)\cup\suppt(y,B)$.

Suppose now that $\delta$ is a limit ordinal. Define $A_\delta=\bigcup_{\xi<\delta}A_\xi$. Clearly, conditions (i) and (ii) are satisfied. Condition (iv) is obvious, so it remains to check (iii). Fix $y\in D$. There is nothing to prove if $\delta$ has a countable cofinality, because then $\suppt(y,A_\delta)=\bigcup_{\xi<\delta}\suppt(y,A_\xi)$. Assume $\cf\delta>\aleph_0$.
Since $A_\delta\subs P_\delta X$, we see that $\suppt(y,A_\delta)=\suppt(P_\delta^*y,A_\delta)$. Thus, we may assume that $y=P_\delta^* y$. Continuity of the projectional sequence implies that $y=\lim_{\xi<\delta}P_\xi^*y$, where the limit is with respect to the \weakstar\ topology. We know that $P_\delta^*D$ is contained in a $\Sig$-space, therefore it is \weakstar\ countably tight. It follows that there exists $\al<\delta$ such that $P_\xi^* y=y$ for $\al\loe \xi <\delta$. In particular, $\suppt(y,A_\delta)=\suppt(y,A_\al)$, because if $a\in A_\delta\setminus A_\al$ then $y(a)=(P_\al^*y)a=y(P_\al a)=0$, by (iv). By the induction hypothesis, $\suppt(y,A_\delta)$ is countable. This shows (iii).

Finally, set $A=\bigcup_{\al<\kappa}A_\al$. By (i), $A$ is linearly dense in $X$. It remains to check that $\suppt(y,A)$ is countable for every $y\in D$. Fix $y\in D$. By the assumption, $y\in P_\al^*D$ for some $\al<\kappa$. Thus $P_\xi^*y=y$ whenever $\al\loe\xi<\kappa$ and we conclude like in the limit case of the above construction. This completes the proof.
\end{pf}

The above theorem should be compared with the results of S. Gul$'$ko (see \cite{Gul77, Gul79, Gul98}), where similar preservation was proved for topological spaces which have a continuous injection into a $\Sig$-product.

\begin{wn}
Let $X$ be a Banach space with a projectional sequence $\sett{P_\al}{\al<\kappa}$ such that $P_\al X$ is weakly Lindel\"of determined for each $\al<\kappa$. Then $X$ is a Plichko space.
\end{wn}

\begin{wn}
Given a Banach space $X$, the following properties are equivalent.
\begin{enumerate}
	\item[(a)] ${X^*}$ generates projections in $X$.
	\item[(b)] $X$ is weakly Lindel\"of determined.
\end{enumerate}
\end{wn}

By a result of Orihuela, Schachermayer and Valdivia \cite{OriSchVal91}, the above properties are also equivalent to ``$\pair{\ubal{X^*}}{\weakstar}$ is Corson compact".

\index{Banach space!-- Asplund}\index{Asplund space}
It is natural to ask when $X$ (as a subspace of $X^{**}$) generates projections in $X^*$.
By a result of Fabi\'an and Godefroy \cite{FabGod88} this is the case when $X$ is Asplund. Specifically, assuming $X$ is an Asplund space, the authors of \cite{FabGod88} construct a projectional generator $\pair X\Phi$ in $X^*$. 
On the other hand, Orihuela and Valdivia noted in \cite[Thm. 3]{OriVal90} that the existence of a projectional generator with domain $X$ and with values in $X^*$ implies that $X$ is Asplund. Recall that a Banach space $X$ is {\em Asplund} if the dual of every separable subspace of $X$ is separable. Assume $X$ generates projections in $X^*$ and fix a separable subspace $Y$ of $X$. Fix a countable $M\rloe H(\theta)$ such that $X\in M$ and $Y\cap M$ is dense in $Y$ and let $\map P{X^*}{X^*}$ be the projection with $\im P = \cl(X^*\cap M)$ and $\ker P=\rightorth{X\cap M}$. Then $P^*y=y$ for every $y\in Y$, because $Y\subs\clstar(X\cap M)$. Fix $\phi\in Y^*$ and let $\psi\in X^*$ be an extension of $\phi$. Then
$(P\psi)y = \psi(P^*y) = \psi(y) = \phi(y)$
%$$\pair y{P\psi} = \pair {P^*y}\psi = \pair y\psi = \pair y\phi,$$
for every $y\in Y$. Thus, $\phi = (P\psi)\rest Y$. It follows that $Y^*$ is separable because so is $\im P$. Summarizing, we have:

\begin{prop}
Given a Banach space $X$, the following properties are equivalent.
\begin{enumerate}
	\item[(a)] $X$ is Asplund.
	\item[(b)] $X$ generates projections in $X^*$.
\end{enumerate}
\end{prop}

%We do not know whether, for an arbitrary Banach space $X$ and a norming space $D\subs X^*$, the existence of a projectional generator $\pair D\Phi$ is implied by the property that $D$ generates projections in $X$.

\subsection{A characterization of Plichko spaces}

\begin{tw}\label{weuiohqwo}
Let $X$ be a Banach space and let $r\goe1$. The following properties are equivalent.
\begin{enumerate}
	\item[(a)] $X$ has a commutative $r$-projectional skeleton.
	\item[(b)] $X$ is an $r$-Plichko space.
\end{enumerate}
\end{tw}

\begin{pf}
Implication (b)$\implies$(a) is contained in Proposition \ref{ofhwoehfr}. For the converse implication, we use Theorem \ref{qwrorqw}, Lemma \ref{iuwbfuiq} and induction on the density of $X$. Suppose we have proved that (a)$\implies$(b) for spaces of density $<\kappa$ and fix a Banach space $X$ of density $\kappa$ with a commutative $r$-projectional skeleton $\sett{P_s}{s\in\Gam}$. By (the proof of) Theorem \ref{qwrorqw}, there exists an $r$-projectional resolution of the identity $\s=\sett{P_\al}{\al<\kappa}$ on $X$ such that for each $\al<\kappa$ there is a directed set $S_\al\subs\Gam$ with $P_\al x=\lim_{s\in S_\al}P_sx$ for $x\in X$ (to be formal, we need to assume that $\Gam\cap\kappa=\emptyset$).
Observe that, by continuity, $P_s \cmp P_\al = P_\al\cmp P_s$ holds for every $s\in\Gam$ and $\al<\kappa$. Let $D$ be the norming space induced by $\s$. Fix $y\in D$ and fix $\al<\kappa$. Let $s\in\Gam$ be such that $y=P_s^*y$. Then $P_\al^* y = P_\al^*P_s^* y = P_s^*P_\al^*y\in D$. Hence $P_\al^*D\subs D$. Now use Theorem \ref{popiooier}. In case where $\cf\kappa=\aleph_0$, it may happen that $D\ne \bigcup_{\al<\kappa}P_\al^*D$, but we may replace $D$ by $\bigcup_{\al<\kappa}P_\al^*D$, still having an $r$-norming space.
By Theorem \ref{popiooier}, $\pair XD$ is a Plichko pair, therefore $X$ is $r$-Plichko.
\end{pf}

\section{Spaces of continuous functions}

We discuss a natural class of compact spaces $K$ for which $\cee K$ has a projectional skeleton.

\index{retractional skeleton}\index{r-skeleton}
Let \Rzero\ denote the class of all compacta which have a {retractional skeleton}. Following \cite{KubMich}, a {\em retractional skeleton} (briefly: {\em r-skeleton}) in a compact space $K$ is a family of retractions $\sett{r_s}{s\in\Gam}$ onto metrizable subsets, indexed by an up-directed poset $\Gam$, satisfying the following conditions:
\begin{enumerate}
	\item $s\loe t \implies r_s = r_s\cmp r_t = r_t\cmp r_s$.
	\item For every $x\in X$, $x = \lim_{s\in \Gam}r_s(x)$.
	\item Given $s_0<s_1<\dots$ in $\Gam$, $t=\sup_{\ntr}s_n$ exists and $r_t(x)=\lim_{n\to\infty}r_{s_n}(x)$ for every $x\in K$.
\end{enumerate}
It has been proved in \cite{KubMich} that Valdivia compacta are precisely those compact spaces which have a commutative r-skeleton. The ordinal $\omega_2+1$ is an example of a space in class \Rzero\ which is not Valdivia compact. It is clear that every r-skeleton induces a $1$-projectional skeleton on the space of continuous functions; that is:

\begin{prop}\label{weoptjperjqr}
Let $K$ be a compact space with a retractional skeleton $\sett{r_s}{s\in\Gam}$. Then $\sett{r_s^*}{s\in\Gam}$ is a projectional skeleton in $\cee K$, where $r_s^*$ denotes the transformation adjoint to $r_s$.
\end{prop}

A simple application of Lemma \ref{iuwbfuiq} shows that every space from class \Rzero\ can be decomposed into a continuous inverse sequence of retractions onto smaller spaces in class \Rzero\ (notion dual to a PRI). This shows that $\Rzero\subs \R$, where \R\ is the smallest class of spaces containing all metric compacta and closed under limits of continuous inverse sequences of retractions (see \cite{BurKubTod, Kub06}). Note that class \Rzero\ restricted to spaces of weight $\loe\aleph_1$ coincides with the class of Valdivia compacta (\cite[Cor. 4.3]{KubMich}). This is not the case with class \R\ (see \cite[Example 4.6(b)]{KubMich}), therefore $\Rzero\ne \R$. 

Let us admit that the converse to Proposition \ref{weoptjperjqr} fails, namely there exist compact spaces $K\notin\R$ such that $\cee K$ is $1$-Plichko, see \cite{BanKub}.
On the other hand, by Lemma \ref{eiurqrad}, we have:

\begin{prop}
Let $D\subs X^*$ be a $1$-norming space which generates projections in a Banach space $X$. Then $\ubal{X^*}$ endowed with the \weakstar\ topology belongs to class \Rzero.
\end{prop}

The following results are dual to Theorems \ref{wejirpqwrt}, \ref{eoghjnoasf} and \ref{rwqroiufd} respectively.

\begin{tw}\label{aefjqjipf} Let $K$ be a compact space and let $D\subs K$ be a dense countably closed set. The following properties are equivalent:
\begin{enumerate}
  \item[(a)] $K\in\Rzero$ and $D$ is induced by an r-skeleton in $K$.
  \item[(b)] For every sufficiently big cardinal $\theta$, for every countable elementary substructure $M$ of $H(\theta)$ with $K,D\in M$, the quotient $\map{q^M_K}K{K\by M}$ restricted to $\cl(D\cap M)$ is one-to-one.
\end{enumerate}
\end{tw}

\begin{pf}
Assume (a) and fix a countable $M\rloe\pair{H(\theta)}\in$ such that $K,D\in M$. By elementarity, there exists $\sett{r_s}{s\in\Gam}\in M$ which is an r-skeleton in $K$ such that $D=\bigcup_{s\in \Gam}\img{r_s}K$. Fix $x,y\in\cl(D\cap M)$, $x\ne y$. Let $t=\sup(\Gam\cap M)$. By elementarity, $D\cap M\subs\img {r_t}K$, therefore also $\cl(D\cap M)\subs\img{r_t}K$. It follows that $x=r_t(x)$ and $y=r_t(y)$. Let $\sett{s_n}{\ntr}\subs \Gam\cap M$ be increasing and such that $t=\sup_{\ntr}s_n$. Then $x=\lim_{n\to\infty}r_{s_n}(x)$ and $y=\lim_{n\to\infty}r_{s_n}(y)$. It follows that $r_{s_k}(x)\ne r_{s_k}(y)$ for all but finitely many $k\in\nat$. Since $r_{s_k}\in M$ and $\img{r_{s_k}}K$ is second countable, we deduce that $x\not\sim_My$.
This shows that (a)$\implies$(b).

The proof of (b)$\implies$(a) is similar to that of (b)$\implies$(a) in Theorem \ref{wejirpqwrt}: the family $\Gam$ of all countable $M\rloe\pair{H(\theta)}\in$ with $K,D\in M$ is an r-skeleton on $K$, since each $q^M_K$ can be treated as a retraction of $K$ onto $\cl(K\cap M)$.
\end{pf}

Note that countably tight spaces in class \Rzero\ are precisely Corson compacta. Thus, in case $K$ is Corson compact, we have $D=K$ and the above theorem gives Bandlow's characterization \cite{Ban91}.

\begin{prop}
Class \Rzero~is closed under arbitrary products.
\end{prop}

\begin{pf} Let $\setof{K_\al}{\al\in\kappa}$ be a family of spaces in \Rzero~and let $K=\prod_{\al\in\kappa}K_\al$. For each $\al\in\kappa$ choose a dense set $D_\al\subs K_\al$ which is induced by a fixed r-skeleton in $K_\al$. Fix $p\in\prod_{\al\in\kappa}D_\al$ and let $D$ be the $\Sig$-product of $\sett{D_\al}{\al\in\kappa}$ based on $p$, i.e.
$$D=\setof{x\in\prod_{\al\in\kappa}D_\al}{|\suppt_p(x)|\loe\aleph_0},$$
where $\suppt_p(x):=\setof{\al\in\kappa}{x(\al)\ne p(\al)}$.
It is clear that $D$ is countably closed and dense in $K$. We check condition (b) of Theorem \ref{aefjqjipf}. Let $\theta>\kappa$ be a regular cardinal such that for every $\al<\kappa$ statement (b) of Theorem \ref{aefjqjipf} holds for every countable $M\rloe H(\theta)$ with $K_\al,D_\al\in M$. Fix a countable $M\rloe H(\theta)$ with $\sett{K_\al}{\al\in\kappa},\sett{D_\al}{\al\in\kappa}\in M$. Let $S=\kappa\cap M$. Observe that
$$\cl(D\cap M) \subs \setof{x\in K}{\suppt_p(x)\subs S}.$$
Fix $x\ne y$ in $\cl(D\cap M)$. Then $x(\al)\ne y(\al)$ for some $\al$. By the above remark, $\al\in M$. Thus $K_\al,D_\al\in M$ and therefore $q^M_{K_\al}$ is one-to-one on $\cl(D_\al\cap M)$. Finally, if $f\in \cee{K_\al}\cap M$ separates $x(\al)$ from $y(\al)$, then $g=f\cmp \pr_\al\in M$ and $g(x)\ne g(y)$, where $\pr_\al$ denotes the projection onto the $\al$-th coordinate. This shows that $q^M_K$ is one-to-one on $\cl(D\cap M)$. By Theorem \ref{aefjqjipf}, $K\in\Rzero$.
\end{pf}

\begin{tw}\label{oituweghsoewt}
Assume $\sett{R_s}{s\in\Gam}$ is an r-skeleton in a compact space $K$ and let $D=\bigcup_{s\in\Gam}\img{R_s}K$. Then
\begin{enumerate}
	\item[$(1)$] $D$ is dense in $K$ and for every countable set $A\subs D$ the closure $\cl_K(A)$ is metrizable and contained in $D$.
	\item[$(2)$] $D$ has a countable tightness.
	\item[$(3)$] $D$ is a normal space and $K=\beta D$.
\end{enumerate}
\end{tw}

\index{r-skeleton!-- dense set induced by}
We shall say that $D$ is the {\em dense set induced by} $\sett{R_s}{s\in\Gam}$.

\begin{pf} (1) follows from the $\sig$-directedness of $\Gam$: every countable subset of $D$ is contained in some $K_s:=\img{R_s}K$. 
The countable tightness of $D$ follows from Theorem \ref{rwqroiufd}(b), because $D\subs {\cee K}^*$ generates projections in $\cee K$.

It remains to prove that $D$ is a normal space and that $K=\beta D$. Fix disjoint relatively closed sets $A,B\subs D$. We claim that $\cl_K(A)\cap \cl_K(B)=\emptyset$. This will also show that $K=\beta D$. Suppose $p\in \cl_K(A)\cap \cl_K(B)$ and fix a countable $M\rloe H(\theta)$ (where $\theta$ is sufficiently big) such that $A,B,p,\sett{R_s}{s\in\Gam}\in M$. Let $\delta=\sup(\Gam\cap M)$. Then $R_\delta(p)\in D$, so we may assume that $R_\delta(p)\notin A$ (interchanging the roles of $A$ and $B$, if necessary). 
Recall that $K_\delta:=\img{R_\delta}K$ is the limit of inverse system $\invsys KR{\Gam\cap M} st$, where $R^t_s=R_s\rest K_t$ (see \cite[Lemma 3.4]{KubMich}). Thus, there are $t\in \Gam\cap M$ and an open set $V\subs K_t$ such that $U:=K_\delta\cap \inv{(R_t)}V$ is a neighborhood of $p$ in $K_\delta$ disjoint from $A$. On the other hand, $\inv{(R_t)}V\cap A\nnempty$. Since $K_t$ is second countable, we may assume that $V\in M$. Thus, by elementarity, there is $a\in M$ such that $a\in \inv{(R_t)}V\cap A$. Finally, $a\in K_\delta$, so $a\in U\cap A$, a contradiction.
\end{pf}

A Banach space analogue of part (3) of the above result looks as follows.

\begin{prop}
Let $D$ be a norming space induced by a projectional skeleton $\sett{P_s}{s\in\Gam}$ in a Banach space $X$. Then for every \weakstar\ continuous function $\map fD\Err$ there exists $t\in\Gam$ such that $f = f\cmp P_t^*\rest D$. 
\end{prop}

\begin{pf}
Fix $n>0$ and consider $K_n=n\ubal{X^*}$. Then $\sett{P_s^*\rest K_n}{s\in\Gam}$ is a retractional skeleton in $K_n$. By \cite[Lemma 5.1]{KubMich}, there exists $s_n\in\Gam$ such that $f_n = f\cmp P_{s_n}^*\rest (D\cap K_n)$. We may assume that $s_1\loe s_2\loe \dots$. Let $t=\sup_{\ntr}s_n$. Then $f = f\cmp P_t^*\rest D$.
\end{pf}

We are now able to determine \weakstar\ compact subsets of spaces induced by projectional skeletons.

\begin{tw}
Assume $D\subs X^*$ generates projections in a Banach space $X$. Then every compact subset of $D$ is Corson.
\end{tw}

\begin{pf}
Let $K\subs D$ be compact with respect to the \weakstar\ topology. We use Bandlow's characterization \cite{Ban91}, which is a special case of Theorem \ref{aefjqjipf}.
Fix a suitable countable $M\rloe\pair{H(\theta)}\in$ and fix $p\ne q$ in $\clstar(K\cap M)$. Then there is $x\in X$ such that $p(x)\ne q(x)$. Since $\pair XD$ has Property $\Omega$ (see Proposition \ref{weriojqwpriq}), there exists $y\in \cl(X\cap M)$ such that $x-y\in\leftorth{D\cap M}$. In particular, $p(x)=p(y)$ and $q(x)=q(y)$. Now, the continuity of $p$ and $q$, find $z\in X\cap M$ such that $p(z)\ne q(z)$. Thus the function $\phi\mapsto\phi(z)$ is an element of $M$ which separates $p$ and $q$. This shows that $p\not\sim_Mq$. Finally, Bandlow's theorem \cite{Ban91} (or a special case of Theorem \ref{aefjqjipf}) shows that $K$ is Corson.
\end{pf}

We do not know whether the converse holds, namely whether every norming \weakstar\ Corson compact set generates projections, see Question \ref{wepfjpqrewer}. 

A preservation theorem for Valdivia compacta, dual to Theorem \ref{popiooier}, looks as follows. We omit its proof, since it can be easily deduced from (the proof of) Theorem \ref{popiooier}.
 
\index{Valdivia pair}
Given a Valdivia compact $K$, let us call $\pair KD$ a {\em Valdivia pair} if $D$ is dense in $K$ and there is an embedding $\map jK{[0,1]^\kappa}$ such that $\img jD\subs\Sig(\kappa)$.

\begin{tw}
Let $\sett{r_\al}{\al<\kappa}$ be a continuous retractive sequence in a compact space $K$. Let $D\subs K$ be a dense set such that for each $\al<\kappa$, $\pair{\img{r_\al}K}{\img{r_\al}D}$ is a Valdivia pair
and $\img{r_\al}D\subs D$ for every $\al<\kappa$. Then $\pair KD$ is a Valdivia pair.
\end{tw}

The above result leads to another proof of \cite[Thm. 6.1]{KubMich}, saying that a compact space with a commutative r-skeleton is Valdivia compact. Another corollary is the following.

\begin{wn}
The limit of a continuous retractive inverse sequence of Corson compacta is Valdivia compact.
\end{wn}

\section{Final remarks and open problems}\label{finalremarks}

As we have already mentioned, the ordinal $\omega_2+1$ provides an example of a compact space in class \Rzero\ whose space of continuous functions is not Plichko. An r-skeleton in $\omega_2+1$ can be constructed as follows. Denote by $\Gam$ the family of all countable closed subsets $A$ of $\omega_2$ such that $0\in A$ and every isolated point of $A$ is isolated in $\omega_2$. Given $A\in \Gam$, define $\map{r_A}{\omega_2+1}{\omega_2+1}$ by setting $r_A(\al)=\max(A\cap[0,\al])$. It is straight to check that $r_A$ is a retraction onto $A$ (continuity follows from the assumption concerning isolated points). It is easy to check that $\ft=\sett{r_A}{A\in\Gam}$ is an r-skeleton. Obviously, this skeleton is not commutative. On the other hand, the dual $\cee{\omega_2+1}^*$ is $1$-Plichko (see \cite[Example 4.10(a)]{Kal01} or \cite[Example 6.10]{Kal00a}).

\begin{ex}\label{Pqrqwrwwegwe}
There exists a $1$-projectional skeleton $\s$ on $\ell_1(\omega_2)$ such that $D(\s)$ is not a $\Sigma$-space, i.e. $\pair{\ell_1(\omega_2)}{D(\s)}$ is not a Plichko pair.
\end{ex}

\begin{pf}
We shall use $\omega_2+1$ instead of $\omega_2$ as the coordinate set, because $\omega_2+1=[0,\omega_2]$ has the maximal element with respect to the natural well order.
Let $\Gam$ consist of all countable subsets $S$ of $\omega_2+1$ such that $\omega_2\in S$. Given $S\in\Gam$, define $\map{f_S}{\omega_2+1}{\omega_2+1}$ by $f_S(\al) = \min(S\cap[\al,\omega_2])$. Note that $f_S$ is generally discontinuous with respect to the interval topology on $\omega_2+1$. Further, define $\map{Q_S}{\ell_1(\omega_2+1)}{\ell_1(\omega_2+1)}$ by setting $(Q_Sx)(\al) = \sum_{\xi\in f_S^{-1}(\al)}x(\xi)$. It is clear that $Q_S$ is a well defined linear projection onto $\ell_1(S)\subs\ell_1(\omega_2+1)$. Further, $\norm{Q_S}=1$. Fix $S\subs T$ in $\Gam$.
Clearly, $Q_T\cmp Q_S = Q_S$, because $Q_T$ is identity on $\ell_1(T)\sups\ell_1(S)$.
Now observe that $f_S\cmp f_T = f_S$. Thus, given $x\in \ell_1(\omega_2+1)$ and $\al\in\omega_2+1$, we have that $f_S^{-1}(\al) = \bigcup_{\xi\in f_S^{-1}(\al)}f_T^{-1}(\xi)$ and consequently
$$(Q_S Q_T x)(\al) = \sum_{\xi\in f_S^{-1}(\al)} (Q_T x)(\xi) = \sum_{\xi\in f_S^{-1}(\al)} \; \sum_{\eta\in f_T^{-1}(\xi)} x(\eta) = \sum_{\xi\in f_S^{-1}(\al)} x(\xi) = (Q_S x)(\al).$$
Finally, given $S_0\subs S_1\subs \dots$ in $\Gam$, the set $S_\infty = \bigcup_{\ntr}S_n$ is an element of $\Gam$ and $\bigcup_{\ntr}\im(Q_{S_n})$ is clearly dense in $\im(Q_{S_\infty})$.
It follows that $\s=\sett{Q_S}{S\in\Gam}$ is a projectional skeleton in $\ell_1(\omega_2+1)$.

Now suppose that $\pair{\ell_1(\omega_2)}{D(\s)}$ is a Plichko pair. By Theorem \ref{weuiohqwo}, there exists a commutative projectional skeleton $\ft = \sett{P_t}{t\in\Delta}$ such that $D(\ft)=D(\s)$. By Lemma \ref{weiopfjpqw}, there exists a cofinal subset $\Gam'\subs\Gam$ such that for every $S\in \Gam'$ there is $t=t(S)\in \Delta$ with $Q_S=P_t$. In particular, $Q_S\cmp Q_T = Q_T\cmp Q_S$ whenever $S,T\in \Gam'$. We shall derive a contradiction by finding $S,T\in \Gam'$ such that $Q_S\cmp Q_T \ne Q_T\cmp Q_S$. 

Given $S\in\Gam'$, define $\phi(S)=\sup(S\cap \omega_2)$. Construct a chain $\sett{S_\al}{\al<\omega_1}$ in $\Gam'$ so that $\phi(\al)<\phi(\beta)$ whenever $\al<\beta$. This is possible, because $\Gam'$ is cofinal in $\Gam$.
Let $\delta=\sup_{\al<\omega_1}\phi(\al)$. Fix $T\in\Gam'$ such that $\phi(T)>\delta$. Then $\sup(T\cap\delta)<\delta$, because $\delta$ has cofinality $\omega_1$. Find $\al<\omega_1$ such that $\sup(T\cap \delta)<\phi(S_\al)<\delta$. Let $S=S_\al$. Fix $\xi\in S_\al$ such that $\sup(T\cap \delta)<\xi$. Then $f_T(\xi)>\delta$ and hence $f_S f_T(\xi)=\omega_2$, because $\phi(S)<\delta$. On the other hand, $f_S(\xi)=\xi<\phi(T)$ and hence $f_T f_S(\xi)\loe\phi(T)<\omega_2$. It follows that $f_S f_T(\xi)\ne f_T f_S(\xi)$. Considering the characteristic function of $\sn\xi$ as an element of $\ell_1(\omega_2+1)$, we conclude that $Q_S\cmp Q_T \ne Q_T\cmp Q_S$.
\end{pf}

Given a norming space $D\subs X^*$, let $\Tau_D$ denote the topology on $X$ induced by $D$, i.e. $\Tau_D = \sig(X,D)$. It can be shown that $\pair X{\Tau_D}$ is Lindel\"of, whenever $D$ generates projections in $X$. On the other hand, by the result of Kalenda \cite{Kal00},
$D$ is a $\Sig$-space iff $\pair X{\Tau_D}$ is primarily Lindel\"of and $D\cap \ubal{X^*}$ is \weakstar\ countably compact.

\begin{problem}
Find a topological property of $\pair X{\Tau_D}$ which says when $D$ generates projections in $X$.
\end{problem}

\begin{question}
Assume $X$ has a $1$-projectional skeleton. Does $X$ have a projectional generator?
\end{question}

\begin{question}
Assume $X$ is a Banach space of density $>\aleph_1$ and $\Ef$ is a directed family of $1$-complemented separable subspaces such that $X=\bigcup\Ef$ and $\cl(\bigcup_{\ntr}F_n) \in \Ef$ whenever $\setof{F_n}{\ntr}\subs \Ef$. Does $X$ necessarily have a projectional skeleton?
\end{question}

Note that if $X$ has density $\aleph_1$ then the above assumptions imply the existence of a PRI, see \cite[Lemma 6.1]{Kub07}.

\begin{question}
Let $X$ be a Banach space with a projectional skeleton. Does every closed subspace of $X$ have the separable complementation property?
\end{question}

Note that, by the main result of \cite{Kub07}, a closed subspace of a Plichko space of density $\aleph_1$ may not have a projectional skeleton.

The following question has already been asked by Ond\v rej Kalenda \cite{Kal02}.

\begin{question}
Is $\cee{\omega_2+1}$ embeddable into a Plichko space?
\end{question}

\begin{question}\label{wepfjpqrewer}
Assume $K\subs X^*$ is Corson compact and norming for $X$. Does $K$ generate projections in $X$?
\end{question}

If $K$ is a Corson compact in the dual of a Banach space $X$ and $K$ is norming for $X$, then $X$ embeds into $\cee K$. Thus, if $\cee K$ is WLD then so is $X$ and consequently $K$ generates projections in $X$. It follows that the above question at least consistently has affirmative answer: it is relatively consistent with the usual axioms of set theory that $\cee K$ is WLD for every Corson compact $K$ (see e.g. \cite[Remark 3.2.3)]{ArgMerNeg}).

\section*{Acknowledgements}

The author would like to thank the referee for many valuable remarks, in particular for providing a simpler proof of Lemma \ref{iuwbfuiq}, suggesting Example \ref{Pqrqwrwwegwe} and for pointing out references \cite{Kal01, SimYos}.

\def\cprime{$'$}

\printindex
\end{document}